\documentclass[11pt, keyword]{article}
\usepackage{graphicx, comment}
\usepackage{times}
\usepackage{amsmath}
\usepackage{amssymb}
\usepackage{latexsym}
\usepackage{amsthm}
\newtheorem{theorem}{Theorem}[section]
\newtheorem{lemma}[theorem]{Lemma}
\newtheorem{proposition}[theorem]{Proposition}
 
\newtheorem{definition}[theorem]{Definition}
\newtheorem{corollary}[theorem]{Corollary}
\newtheorem{remark}[theorem]{Remark}
\newcommand{\As}{\mathcal{A}_{\mathrm{s}}}
\newcommand{\Aw}{\mathcal{A}_{\mathrm{w}}}

\newcommand{\ds}{\mathrm{d}_{\mathrm{s}}}
\newcommand{\dw}{\mathrm{d}_{\mathrm{w}}}

\newcommand{\Dc}{\mathcal{E}}
\newcommand{\s}{\mathrm{s}}
\newcommand{\w}{\mathrm{w}}
\newcommand{\Xw}{X_{\mathrm{w}}}
\newcommand{\Xs}{X_{\mathrm{s}}}
\newcommand{\Hw}{H_{\mathrm{w}}}
\newcommand{\Ab}{\mathcal{A}_{\bullet}}
\newcommand{\db}{\mathrm{d}_{\bullet}}
\newcommand{\dd}{\mathrm{d}}
\newcommand{\I}{\mathcal{I}}
\newcommand{\Bb}{B_{\bullet}}
\newcommand{\Bw}{B_{\mathrm{w}}}
\newcommand{\wb}{{\omega}_{\bullet}}
\newcommand{\ww}{{\omega}_{\mathrm{w}}}
\newcommand{\ws}{{\omega}_{\mathrm{s}}}

\newcommand{\ddt}{\frac{d}{dt}}
\newcommand{\Kp}{\mathcal{K}^+}
\newcommand{\K}{\mathcal{K}}
\newcommand{\F}{\mathcal{F}^+}
\newcommand{\Fc}{\mathcal{F}}
\newcommand{\At}{\mathfrak{A}}
\newcommand{\Rc}{\widetilde{R}}

\def\setmarsing{
\oddsidemargin-0in
\evensidemargin-0in
\textwidth5.5in
\textheight8in
}
\setmarsing

\title{Global attractors of evolutionary systems}
\author{Alexey Cheskidov\\
\small
Department of Mathematics, University of Michigan,
Ann Arbor, MI 48109}
\begin{document}
\maketitle
\abstract{
An abstract framework for studying the asymptotic behavior
of a dissipative evolutionary system $\Dc$ with respect to
weak and strong topologies was introduced in \cite{CF} primarily to
study the long-time behavior of the 3D Navier-Stokes equations (NSE) for which
the existence of a semigroup of solution operators is not known.
Each evolutionary system possesses a global attractor in the weak topology,
but does not necessarily in the strong topology. 
In this paper we study the structure of a global attractor for an abstract evolutionary system,
focusing on omega-limits and attracting, invariant, and quasi-invariant sets.
We obtain weak and strong uniform tracking properties of omega-limits and global
attractors. In addition, we discuss
a trajectory attractor for an evolutionary system and
derive a condition under which the convergence to
the trajectory attractor is strong.
}

\section{Introduction}
Existence of a global attractor is a significant feature of many dissipative
partial differential equations (PDEs). 
Often the evolution of solutions to a dissipative PDE
can be described by a semigroup of solution operators.
If the semigroup is asymptotically compact, then the classical theory of semiflows
yields the existence of a compact global attractor (see Hale \cite{H}, Ladyzhenskaya \cite{L}, or Temam \cite{T2}). However, for some PDEs the semigroup
is not asymptotically compact. Some PDEs, such as the 3D Navier-Stokes equations (NSE),
do not even possess the semigroup due to the lack of uniqueness (or proof of uniqueness).

There are several abstract frameworks for studying dynamical systems without
uniqueness.
See Caraballo, P. Mar'n-Rubio \& Robinson \cite{CMR} for a comparison of two canonical ones by
Melnik \& Valero \cite{MV} and Ball \cite{B1}. In the first approach, used by
Babin \& Vishik \cite{BV}, and which goes all the way back to the work by Barbashin \cite{Bar},
a trajectory is a function of time with values in the set of all subsets of a phase space.
The evolutionary system $\Dc$
considered in this paper 
is closer to the Ball's generalized semiflow $G$, where a trajectory is a function
of time with values in the phase space and there may be more than one
trajectory with given initial data.
Since $\Dc$ does not contain the hypotheses of
concatenation and upper semicontinuity with respect to initial data,
the Leray-Hopf weak solutions of the 3D NSE form an evolutionary
system.


In fact, the notion of an evolutionary system, introduced in Cheskidov \& Foias \cite{CF},
 was motivated by  the 3D NSE,
which possesses a global attractor with respect to the weak topology of the natural phase space. This weak global attractor, introduced by Foias and Temam
in \cite{FT85}, captures the
long-time behavior of all Leray-Hopf weak solutions. In particular,  it includes
the support of any time-average measure of the 3D NSE (see
Foias, Manley, Rosa \& Temam \cite{FMRT}).




A phase space $X$ for the evolutionary system $\Dc$ is a metric space (whose metric is
called strong), which is compact in some weaker metric. For a dissipative PDE,
the space $X$ is defined to be an absorbing ball, and the weak metric is a metric
induced by the weak topology. A global attractor for $\Dc$ is the minimal closed
attracting set in the corresponding topology. In \cite{CF} it is shown
that the weak global attractor always exists, it is the maximal invariant set, and 
if the strong global attractor exists, then its weak closure is the weak global attractor.
Moreover, if $\Dc$ is asymptotically compact, then
the weak global attractor becomes the strong compact global attractor.
Applied to the 3D NSE, this result implies the existence of a strong compact global attractor in the case where solutions on the weak global attractor are
continuous in $L^2$ (see Ball \cite{B1} and Rosa \cite{R} for similar results).

In this paper we continue investigating properties of the evolutionary system $\Dc$, 
concentrating on omega-limits and attracting, invariant, and quasi-invariant sets.
Assume that the evolutionary system is not asymptotically compact. Since the phase
space $X$ is weakly compact,
we will see that the omega-limits and attracting sets with respect to the
weak metric possess familiar properties known from the
classical theory of semiflows. 
Our goal will be to examine the corresponding objects with respect to the strong
metric.

The structure of the paper is as follows. In Section 2 we define the evolutionary
system $\Dc$ and compare it with a semiflow. Attracting sets, global
attractors, and omega-limits are defined and studied in Section 3. We
show the existence of a weak global attractor $\Aw$ and deduce a new necessary
and sufficient condition (in terms of the omega-limits) for the existence of the
strong global attractor $\As$.

In Section 4 we further study the evolutionary system under the condition
that $\Dc$ is asymptotically compact. In this case the situation is the same
as in the classical theory of global attractors. We see that the weak
omega-limit of a set uniformly strongly attracts that set. Moreover, the weak omega-limit coincides
with the strong omega-limit. In particular, we recover a result from \cite{CF},
which says that $\As$ exists, is strongly compact, and coincides with $\Aw$.

Section 5 mostly focuses on invariant and quasi-invariant sets. We study an
evolutionary system with an assumption that the family of all the trajectories
is compact in $C([0,\infty);\Xw)$. For instance, Leray-Hopf
weak solutions of the 3D NSE satisfy this property. We prove that the weak
omega-limit of any set is quasi-invariant (i.e., consists of complete
orbits) and possesses a weak uniform tracking property (Theorem~\ref{t:longtimelim}).
Moreover, in the case where the evolutionary system is asymptotically compact,
omega limits possess a strong uniform tracking property (Theorem~\ref{t:strongtrack}), which holds,
for instance, for the  2D NSE supplemented with appropriate boundary conditions.
This generalizes the tracking property by Langa and Robinson \cite{LR}.
Finally, we prove that $\Aw$ is the strong omega-limit of the phase space $X$,
and if $\As$ exists, then it has to coincide with $\Aw$, i.e., $\As$ is automatically weakly closed.

In Section 6 we make one more step towards abstracting results known
for the 3D NSE. We consider
an evolutionary system that satisfies the energy inequality and the strong
convergence almost everywhere
for a weakly convergent sequence of trajectories. Again,  Leray-Hopf weak
solutions of the 3D NSE satisfy these properties. In this case, the strong
continuity of a trajectory $u(t)$ automatically implies that any sequence of
trajectories that uniformly weakly converges to $u(t)$ converges in fact strongly.
In particular, this yields that $\As$ exists and coincides with $\Aw$ provided
that all the complete trajectories are strongly continuous. See Cheskidov,
Friedlander \& Pavlovi\'c \cite{CFP1} for
an implication of this theorem, which is the existence of a strong global attractor for the
inviscid dyadic model of fluid equations introduced by Friedlander \& Pavlovi\'c \cite{FP}
and Katz \& Pavlovi\'c \cite{KP}. This fact is a result of an anomalous dissipation
due to the loss of regularity of solutions, as it was conjectured by 
Onsager \cite{O} for the 3D Euler equation.
See also Cheskidov \cite{C} for an application of this theorem to a viscous dyadic model.

In Section 7 we study a trajectory attractor $\At$ for an evolutionary system.
The trajectory attractor, a global attractor in the trajectory space, was first
introduced by Sell $\cite{S}$ for the 3D NSE, and further studied in Chepyzhov \& Vishik
\cite {CV} and  Sell \& You \cite{SY}. We show that every evolutionary system possesses a trajectory attractor and discuss its connection to the weak global attractor. Moreover, we prove a strong convergence of the trajectories to $\At$
in the case where all the complete trajectories are strongly continuous.

Finally, in Section 8 we show that the Leray-Hopf weak solutions of the 3D NSE
form an evolutionary system satisfying the above mentioned three additional
properties. Therefore, all the results obtained in this paper
apply to the 3D NSE. 

%


\section{Evolutionary system}
Let $(X,\ds(\cdot,\cdot))$ be a metric space endowed with
a metric $\ds$, which will be referred to as a strong metric.
Let $\dw(\cdot, \cdot)$ be another metric on $X$ satisfying
the following conditions:
\begin{enumerate}
\item $X$ is $\dw$-compact.
\item If $\ds(u_n, v_n) \to 0$ as $n \to \infty$ for some
$u_n, v_n \in X$, then $\dw(u_n, v_n) \to 0$ as $n \to \infty$.
\end{enumerate}
Due to the property 2, $\dw(\cdot,\cdot)$ will be referred to as a weak metric on $X$. Denote by $\overline{A}^{\bullet}$ the closure of a set $A\subset X$
in the topology generated by $\db$.
Note that any strongly compact ($\ds$-compact) set is weakly compact
($\dw$-compact), and any weakly closed set is strongly closed.

Let $C([a, b];X_\bullet)$, where $\bullet = \mathrm{s}$ or $\mathrm{w}$, be the space of $\db$-continuous $X$-valued
functions on $[a, b]$ 
endowed with the metric
\[
\dd_{C([a, b];X_\bullet)}(u,v) := \sup_{t\in[a,b]}\db(u(t),v(t)). 
\]
Let also $C([a, \infty);X_\bullet)$ be the space of $\db$-continuous
$X$-valued functions on $[a, \infty)$
endowed with the metric
\[
\dd_{C([a, \infty);X_\bullet)}(u,v) := \sum_{T\in \mathbb{N}} \frac{1}{2^T} \frac{\sup\{\db(u(t),v(t)):a\leq t\leq a+T\}}
{1+\sup\{\db(u(t),v(t)):a\leq t\leq a+T\}}.
\]
Assume that $u, u_n \in C([T,\infty); X_\bullet)$, $n\in \mathbb{N}$, are such that
$u_n|_{[T_1,T_2]} \to u|_{[T_1,T_2]}$ in $C([T_1,T_2];X_\bullet)$  as $n\to \infty$,
for some $T\leq T_1 \leq T_2$. To simplify the notation, in such cases we will
usually write $u_n \to u$ in $C([T_1,T_2];X_\bullet)$.

To define an evolutionary system, first let
\[
\mathcal{T} := \{ I: \ I=[T,\infty) \subset \mathbb{R}, \mbox{ or } 
I=(-\infty, \infty) \},
\]
and for each $I \subset \mathcal{T}$, let $\mathcal{F}(I)$ denote
the set of all $X$-valued functions on $I$.
\begin{definition} \label{Dc}
A map $\Dc$ that associates to each $I\in \mathcal{T}$ a subset
$\Dc(I) \subset \mathcal{F}$ will be called an evolutionary system if
the following conditions are satisfied:
\begin{enumerate}
\item $\Dc([0,\infty)) \ne \emptyset$.
\item
$\Dc(I+s)=\{u(\cdot): \ u(\cdot -s) \in \Dc(I) \}$ for
all $s \in \mathbb{R}$.
\item $\{u(\cdot)|_{I_2} : u(\cdot) \in \Dc(I_1)\}
\subset \Dc(I_2)$ for all
pairs $I_1,I_2 \in \mathcal{T}$, such that $I_2 \subset I_1$.
\item
$\Dc((-\infty , \infty)) = \{u(\cdot) : \ u(\cdot)|_{[T,\infty)}
\in \Dc([T, \infty)) \ \forall T \in \mathbb{R} \}.$
\end{enumerate}
\end{definition}
We will refer to $\Dc(I)$ as the set of all trajectories
on the time interval $I$. Trajectories in $\Dc((-\infty,\infty))$ will be called complete.
Let $P(X)$ be the set of all subsets of $X$.
For every $t \geq 0$, define a map
\begin{eqnarray*}
&R(t):P(X) \to P(X),&\\
&R(t)A := \{u(t): u\in A, u \in \Dc([0,\infty))\}, \qquad
A \subset X.&
\end{eqnarray*}
Note that the assumptions on $\Dc$ imply that $R(t)$ enjoys
the following property:
\begin{equation} \label{eq:propR(T)}
R(t+s)A \subset R(t)R(s)A, \qquad A \subset X,\quad t,s \geq 0.
\end{equation}

We will also study evolutionary systems $\Dc$ satisfying the following
assumptions:
\begin{itemize}
\item[A1] $\Dc([0,\infty))$ is a compact set in $C([0,\infty); \Xw)$.
\item[A2] (Energy inequality) Assume that $X$ is a bounded set in some uniformly
convex Banach space $H$ with
the norm denoted by $|\cdot|$, such that $\ds(x,y)=|x-y|$ for $x,y \in X$. Assume
also that for any $\epsilon >0$, there exists $\delta$, such that
for every $u \in \Dc([0,\infty))$ and $t>0$,
\[
|u(t)| \leq |u(t_0)| + \epsilon,
\]
for $t_0$ a.e. in $ (t-\delta, t)$.
\item[A3] (Strong convergence a.e.) Let $u,u_n \in \Dc([0,\infty))$, be such that
$u_n \to u$ in $C([0, T];\Xw)$ for some $T>0$. Then
$u_n(t) \to u(t)$ strongly a.e. in $[0,T]$.
\end{itemize}
We will see that all the assumptions A1 -- A3 hold for an evolutionary system
consisting of the Leray-Hopf weak solutions of the 3D Navier-Stokes equations.

Let us now show that a semiflow defines an evolutionary system.
In most applications, the phase space $H$ (a functional space in which trajectories are
defined) is a separable reflexive Banach space.
Consider a semigroup of continuous operators $S(t): H \to H$, $t \geq 0$ satisfying the following properties:
\begin{equation} \label{eq:semigroup}
S(t+s) = S(t)S(s), \quad t,s \geq 0, \qquad S(0)= \mbox{Identity operator.}
\end{equation}
A trajectory $u(t)$ is a mapping from $\mathbb{R^+}$ to $H$, such that
\[
u(t+s) = S(t)u(s), \qquad t,s \geq 0.
\]
A ball $B\subset H$ is called an absorbing ball if for any bounded
set $A \subset H$, there exists $t_0$, such that
\[
S(t)A \subset B, \qquad \forall t \geq t_0.
\]
Assume that the semiflow is dissipative, i.e., there exists an absorbing ball.
Then, if we are interested in a long-time behavior of solutions,
it is enough to consider a restriction of the semiflow to the absorbing ball.
So, we let $X$ be a closed absorbing ball. Since $H$ is a separable reflexive
Banach space, both the strong and the weak topologies on $X$ are metrizable.
Now define the map $\Dc$ in the following way:
\[
\Dc(I) := \{u(\cdot): u(t+s)=S(t)u(s) \mbox{ and }
u(s) \in X \ \forall t \geq 0 , s\in I \}.
\]
Conditions 1--4 in the definition of the evolutionary system $\Dc$ follow
from the semigroup properties (\ref{eq:semigroup}) of $S(t)$.
In addition, let $T$ be such that
\[
S(t)X \subset X \qquad \forall t \geq T.
\]
Then we have
\[
R(t)A = S(t)A, \qquad \forall A \subset S(T)X, \ t \geq 0.
\]

\section{Attracting sets, $\omega$-limits, and global attractors}

For a set $A \subset X$ and $r>0$, denote
$
B_{\bullet}(A,r) = \{u: \ \db(u, A) < r\},
$
where
\[
\db(u, A):=\inf_{x\in A}\db(u,x), \qquad \bullet = \mathrm{s,w}.
\]
A set $A \subset X$ uniformly attracts a set $B \subset X$
in $\db$-metric ($\bullet = \mathrm{s,w}$) if for any $\epsilon>0$
there exists $t_0$, such that
\[
R(t)B \subset B_{\bullet}(A, \epsilon), \qquad \forall t \geq t_0.
\]

\begin{definition}
A set $A \subset X$ is a $\mathrm{d}_{\bullet}$-attracting set
($\bullet = \mathrm{s,w}$) if it uniformly
attracts $X$ in $\mathrm{d}_{\bullet}$-metric.
\end{definition}

\begin{definition}
A set
$\mathcal{A}_{\bullet}\subset X$ is a
$\mathrm{d}_{\bullet}$-global attractor ($\bullet = \mathrm{s,w}$) if
$\mathcal{A}_{\bullet}$ is a minimal $\mathrm{d}_{\bullet}$-closed
$\mathrm{d}_{\bullet}$-attracting  set.
\end{definition}

Note that since $X$ may not be strongly compact, the intersection of two 
strongly closed strongly attracting sets may not be strongly attracting.
Nevertheless, later we will see that if $\Ab$ exists, then it is unique.
%

\begin{definition} The $\wb$-limit ($\bullet= \mathrm{s, w}$)
of a set $A \subset X$ is
\[
\wb (A):=\bigcap_{T\geq 0} \overline{\bigcup_{t \geq T} R(t) A}^{\bullet}.
\]
\end{definition}
An equivalent definition of the $\wb$-limit set is given by
\[
\begin{split}
\wb (A)=\{&x \in X: \mbox{ there exist a sequence } t_n \to \infty \mbox{ as }
 n \to \infty \mbox{ and } x_n \in R(t_n)A, \\
& \mbox{such that } x_n \to x \mbox{ in }
\db\mbox{-metric} \mbox{ as } t_n \to \infty\}.
\end{split}
\]

The following are some properties of $\omega$-limits that immediately
follow from the definition.

\begin{lemma} \label{l:wsww}
Let $A \subset X$. Then
\begin{itemize}
\item[(a)] $\wb(A)$ is $\db$-closed
($\bullet= \mathrm{s, w}$).
\item[(b)] $\ws(A) \subset \ww(A)$.
\item[(c)]  If $\ww(A)$ is strongly compact and uniformly strongly attracts
$A$, then $\ws(A) = \ww(A)$.
\end{itemize}
\end{lemma}
\begin{proof}
Clearly, part (a) follows from the definition. To show part (b), 
take any $x \in \ws(A)$. By the definition of $\ws$-limit, there exist a
sequence $t_n \to \infty$ as $n \to \infty$ and a sequence $x_n \in R(t_n)A$,
such that $x_n \to x$ strongly as $n \to \infty$. In particular, $x_n \to x$ weakly
as $n \to \infty$. Hence, $x \in \ww(A)$, which proves (b).

Now assume that $\ww(A)$ is strongly compact and uniformly strongly attracts
$A$. Take any $x \in \ww(A)$. By the definition of $\ww$-limit, there exist a
sequence $t_n \to \infty$ as $n \to \infty$ and a sequence $x_n \in R(t_n)A$,
such that $x_n \to x$ weakly as $n \to \infty$.
Since $\ww(A)$ strongly attracts $A$,  there exists a sequence
$a_n \in \ww(A)$, such that
\[
\ds(x_n, a_n) \to 0, \qquad \mbox{as} \qquad n \to \infty.
\]
Note that $a_n \to x$ weakly as $n \to \infty$.
Since $\ww(A)$ is strongly compact, this convergence is in fact strong.
Hence, $x_n \to x$ strongly as $n \to \infty$. Therefore, $x \in \ws(A)$,
which proves (c).
\end{proof}

\begin{lemma} \label{l:dbsubseta}
Let $A$ be a $\db$-closed $\db$-attracting set ($\bullet= \mathrm{s, w}$). Then
\[
\wb(X) \subset A.
\]
\end{lemma}
\begin{proof}
Assume that there exists $x \in \wb(X) \setminus A$. Since $A$ is $\db$-closed,
there exists $\epsilon >0$, such that
\begin{equation} \label{eq:temp89}
x \notin B_{\bullet}(A, \epsilon).
\end{equation}
On the other hand, by the definition of $\omega$-limit, there exist a
sequence $t_n \to \infty$ as $n \to \infty$ and
a sequence $x_n \in R(t_n)X$, such that $x_n \to x$ in $\db$-metric
as $n \to \infty$. Since $A$ is $\db$-attracting,
\[
R(t_n)X \subset B_{\bullet}(A, \epsilon/2), 
\]
for $n$ large enough. Therefore, $x_n \in B_{\bullet}(A, \epsilon/2)$ and
consequently $x \in B_{\bullet}(A, \epsilon)$,
which contradicts \eqref{eq:temp89}.
\end{proof}

Now we can show the uniqueness of a global attractor.
\begin{theorem} \label{l:wb=ab}
If $\Ab$ exists ($\bullet= \mathrm{s, w}$), then
\[
\Ab = \wb(X).
\]
\end{theorem}
\begin{proof}
Thanks to Lemma~\ref{l:dbsubseta}, $\wb(X) \subset \Ab$. Assume that
there exists $a\in\Ab \setminus \wb(X)$. Since $a \notin \wb(X)$, there exist $\epsilon>0$ and a time $t_0>0$, such that
\[
R(t)X \cap \Bb(a,\epsilon)= \emptyset, \qquad \forall t>t_0.
\]
Hence, $\Ab \setminus \Bb(a,\epsilon)$ is a $\db$-closed $\db$-attracting set
strictly included in  $\Ab$, which contradicts the definition of a global
attractor.
\end{proof}

The above results imply the following
characterization of the existence of the $\db$-global attractor.

\begin{theorem} \label{c:Aiff}
$\Ab$ exists if and only if $\wb(X)$ is a $\db$-attracting set.
\end{theorem}
\begin{proof}
If $\Ab$ exists, then Theorem~\ref{l:wb=ab} implies that $\wb(X)=\Ab$.
Therefore, $\wb(X)$ is a $\db$-attracting set.

Assume now that $\wb(X)$ is a $\db$-attracting set.
Note that $\wb(X)$ is also $\db$-closed due to Lemma~\ref{l:wsww}.
Then Lemma~\ref{l:dbsubseta} implies that $\wb(X)$ is the minimal
$\db$-closed $\db$-attracting set, i.e., $\wb(X)$ is the $\db$-global attractor.
\end{proof}

%

Now we will study $\ww$-limit sets. The following theorem is an extension
of a well known result for semiflows (see \cite{H,L,T2}).

\begin{theorem} \label{l:oncompactsemigroup}
Let $A \subset X$ be such that there exists $u\in\Dc([0,\infty))$ with $u(0) \in A$. Then $\ww(A)$ is a nonempty weakly compact set.
In addition, $\ww(A)$ uniformly weakly attracts $A$.
\end{theorem}

\begin{proof}
Since $X$ is weakly compact,
\[
W(T):=\overline{\bigcup_{t \geq T} R(t) A}^\mathrm{w}
\]
is a nonempty weakly compact set for all $T\geq 0$. In addition,
$W(s) \subset W(t)$ for all $s\geq t \geq 0$. Thus,
\[
\ww(A)= \bigcap_{T\geq 0} W(T)
\]
is a nonempty  weakly compact set.

We will now prove that $\ww (A)$ uniformly weakly attracts $A$. Assume it
does not. Then there exists $\epsilon>0$, such that
\[
V(t):=W(t)\cap(X\setminus \Bb(\ww(A), \epsilon)) \ne \emptyset, \qquad \forall t\geq 0.
\]
Since $V(t)$ is weakly compact and $V(s) \subset V(t)$ for
$s\geq t \geq 0$, there exists
\[
x \in \bigcap_{t\geq0} V(t).
\]
Hence, $x \in \ww(A)$. However, this together with the definition of $V(t)$
implies that $x\notin V(t)$, $t\geq 0$, a contradiction.
\end{proof}

Finally, from the results of this section we immediately recover the following theorem
from \cite{CF}:

\begin{theorem} \label{thm:exofA}
A weak global attractor $\Aw$ exists. Moreover, if $\As$ exists, then
$\overline{\As}^{\w}=\Aw$.

\end{theorem}
\begin{proof}
Thanks to Theorem~\ref{l:oncompactsemigroup}, $\ww(X)$ is a weakly
closed weakly attracting set. Therefore, $\Aw$ exists and $\Aw=\ww(X)$
due to Theorems~\ref{c:Aiff} and \ref{l:wb=ab}.

Assume now that $\As$ exists.
Since $\overline{\As}^{\w}$ is a strongly attracting set, it is also weakly attracting.
Moreover, since it is weakly closed, Lemma~\ref{l:dbsubseta} implies that
$\ww(X) \subset \overline{\As}^{\w}$. On the other hand, thanks to
Theorem~\ref{l:wb=ab},
$\As=\ws(X)$. Hence, $\ww(X) \subset \overline{\ws(X)}^{\w}$.
Therefore, $\ww(X)=\overline{\ws(X)}^{\w}=\overline{\As}^{\w}$ due to Lemma~\ref{l:wsww}.

\end{proof}

\section{Existence of a strong global attractor}
\begin{definition}
The evolutionary system $\Dc$ is asymptotically compact if for any
$t_n \to \infty$ as $n\to \infty$ and any $x_n \in R(t_n) X$, the  sequence
$\{x_n\}$ is relatively strongly compact.
\end{definition}

\begin{theorem} \label{t:acompomega}
Let $\Dc$ be asymptotically compact.
Let $A \subset X$ be such that there exists $u\in\Dc([0,\infty))$ with $u(0) \in A$. 
Then $\ws(A)$ is a nonempty strongly compact set that uniformly strongly
attracts $A$, and $\ws(A) = \ww(A)$.
\end{theorem}
\begin{proof}
Since there exists $u\in\Dc([0,\infty))$ with $u(0) \in A$,
Theorem~\ref{l:oncompactsemigroup} implies that $\ww(A)$ is nonempty.
First we will show that $\ww(A)$ uniformly strongly attracts $A$. 
Assume that it does not. Then there
exist $\epsilon>0$, $x_n \in X$, and $t_n \to \infty$ as
$n\to \infty$, such that
\begin{equation} \label{e:xnrtnABs}
x_n \in R(t_n)A \setminus B_{\s}(\ww(A), \epsilon), \qquad
\forall n \in \mathbb{N}.
\end{equation}
Since  $\Dc$ is asymptotically compact, we have that
$\{x_n\}$ is relatively strongly compact.
Passing to a subsequence and dropping a subindex, we may assume that
there exists $x\in X$, such that
\begin{equation} \label{eq:dwxnx}
\ds(x_n, x) \to 0 \qquad \mbox{as} \qquad n \to \infty.
\end{equation}
Then also $x_n \to x$ weakly as $n \to \infty$. Therefore,
we have that $x  \in \omega_{\mathrm{w}}(A)$.
Hence, thanks to \eqref{eq:dwxnx},
there exists $n \in \mathbb{N}$, such that
\[
x_n \in B_{\s}(\ww(A), \epsilon),
\]
a contradiction with \eqref{e:xnrtnABs}.

Now note that $\ws(A) \subset \ww(A)$ due to Lemma~\ref{l:wsww},.
On the other hand, let $x \in \ww(A)$. By
the definition of $\ww$-limit, there exist $t_n \to \infty$ as $n \to \infty$ and
$x_n \in R(t_n)A$, such that 
\[
\dw(x_n, x) \to 0 \qquad \mbox{ as } \qquad n \to \infty.
\]
Thanks to the asymptotic compactness of $\Dc$, this convergence
is in fact strong. Therefore, $x\in \ws(A)$. Hence, $\ws(A)= \ww(A)$.

Finally, we have to show that $\ws(A)$ is strongly compact. Take any sequence
$a_n \in \ws(A)$. By the definition of $\ws$-limit,
there exist $t_n \to \infty$ and $x_n \in R(t_n)A$, such that 
\[
\ds(x_n, a_n) \to 0 \qquad \mbox{ as } \qquad n \to \infty.
\]
Note that $\{x_n\}$ is relatively strongly compact due to the asymptotic
compactness
of $\Dc$. Hence, $\{a_n\}$ is relatively strongly compact and, consequently,
$\ws(A)$ is relatively strongly compact. Due to Lemma~\ref{l:wsww}, $\ws(A)$ is
also strongly closed. Therefore, $\ws(A)$ is strongly compact, which concludes
the proof.
\end{proof}

In particular, we automatically have the following result proved in \cite{CF}
for evolutionary systems, which generalizes corresponding results for generalized semiflows
and semiflows \cite{B1, H,HLS,L}.

\begin{theorem} \label{t:asymptoticcompact}
If the evolutionary system $\Dc$ is asymptotically compact, then
$\Aw$ is a strongly compact strong global attractor.
\end{theorem}
\begin{proof}
Since $\Dc([0,\infty)) \ne \emptyset$, thanks to Theorem~\ref{t:acompomega},
$\ws(X)$ is a $\ds$-compact $\ds$-attracting set and $\ws(X)=\ww(X)=\Aw$.
Now Theorems~\ref{c:Aiff} and \ref{l:wb=ab} imply that $\As$ exists and
$\As=\ws(X)=\Aw,$
which concludes the proof.
\end{proof}

\section{Invariance and tracking properties}

In this section we will further study an evolutionary system $\Dc$
satisfying property A1:
\[
\Dc([0,\infty)) \mbox{ is a compact set in } C([0,\infty); \Xw).
\]
In order to extend the notion of invariance from a semiflow to an
evolutionary system, we will need the following mapping:
\[
\Rc(t) A := \{u(t): u(0) \in A, u \in \Dc((-\infty,\infty))\}, \qquad A \subset X,
\ t \in \mathbb{R}.
\]
\begin{definition} A set $A\subset X$ is positively invariant if
\[
\Rc(t) A\subset A, \qquad \forall t\geq 0.
\]
$A$ is invariant if
\[
\Rc(t) A= A, \qquad \forall t\geq 0.
\]
$A$ is quasi-invariant if for every $a\in A$ there exists a complete trajectory
$u \in \Dc((-\infty,\infty))$ with $u(0)=a$ and $u(t) \in A$ for all $t\in \mathbb{R}$.
\end{definition}
Note that the definition of invariance coincides with the classical one
in the case where $\Dc$ is a semiflow or Ball's generalized semiflow.

If $A$ is invariant, then clearly $A$ is quasi-invariant.
Note also that if $A$ is quasi-invariant, then
\begin{equation} \label{e:inclus1}
A \subset \Rc(t) A \subset R(t) A, \qquad \forall t\geq 0.
\end{equation}
This together with Lemma~\ref{l:wsww} imply that 
\begin{equation} \label{e:inclus}
A \subset \ws(A) \subset \ww(A), \qquad A \mbox { is quasi-invariant.}
\end{equation}
\begin{theorem} \label{t:ws=ww,inv}
Let $\Dc$ be an evolutionary system satisfying $A1$.
Then $\ww(A)$ is quasi-invariant for every $A \subset X$.
\end{theorem}
\begin{proof}
Take any $x \in \ww(A)$.
There exist $t_n \to \infty$ as $n\to \infty$ and $x_n \in R(t_n) A$,
such that $x_n \to x$ weakly as $n\to \infty$. Then there exist
$u_n \in \Dc([-t_n, \infty))$ with $u_n(-t_n) \in A$ and
$u_n(0) =x_n$. Because of A1 and the definition of $\Dc$, we have that
$\Dc([-t_n,\infty))$ is compact in $C([-t_n, \infty);\Xw)$ and
\[
\{u|_{[-t_1,\infty)} : u \in \Dc([-t_{n},\infty))\}
\subset \Dc([-t_1,\infty))
\]
for every $n$.
Therefore, passing to a subsequence and dropping a subindex, we obtain
that there exists $u^1\in \Dc([-t_1, \infty))$, such that
\[
u_n|_{[-t_1,\infty)} \to u^1 \qquad \mbox{in}  \qquad
C([-t_1, \infty);\Xw),
\]
as $n\to \infty$. Again passing to a subsequence and dropping a subindex,
we obtain that there exists $u^2\in \Dc([-t_2, \infty))$, such that
\[
u_n|_{[-t_2,\infty)} \to u^2  \qquad \mbox{in} \qquad
C([-t_2, \infty);\Xw),
\]
as $n\to \infty$. Note that $u^1(t)=u^2(t)$ on $[-t_1,\infty)$.
By a standard diagonalization process we infer that there exist
a subsequence of $u_n$, still denoted by $u_n$, and
$u\in \Fc((-\infty,\infty))$, such that $u|_{[-T,\infty)} \in \Dc([-T, \infty))$ and 
$u_n \to u$ in $C([-T, \infty);\Xw)$ as $n\to \infty$ for all $T>0$. 
Note that $u(0)=x$. In addition,
by the definition of $\Dc$ we have that $u \in \Dc((-\infty, \infty))$.

Now take any $t_0 \in \mathbb{R}$.
Note that $u_n(t_0) \to u(t_0)$ weakly as $n\to \infty$. Since $u_n(-t_n) \in A$,
we have that $u_n(t_0) \in R(t_0+t_n)A$ for $t_n \geq -t_0$.
Hence, $u(t_0) \in \ww(A)$, i.e., the complete trajectory $u(t)$ stays on $\ww(A)$
for all time. Therefore, $\ww(X)$ is quasi-invariant.
\end{proof}

Applied to a weak global attractor, this theorem will have several important consequences. 
The following is one of them.

\begin{theorem} \label{c:wwinsideA}
Let $\Dc$ be an evolutionary system satisfying A1. Let $A\subset X$ be such that
$\ww(A) \subset A$. Then $\ws(A)=\ww(A)$.
\end{theorem}
\begin{proof}
Thanks to Theorem~\ref{t:ws=ww,inv}, we have that $\ww(A)$ is quasi-invariant. 
Then \eqref{e:inclus1} implies that $\ww(A) \subset R(t)\ww(A)$  for all $t\geq0$.
Since $\ww(A)\subset A$, it follows that $\ww(A) \subset R(t) A$  for all $t\geq0$.
Therefore, $\ww(A) \subset \ws(A)$. On the other hand,
$\ws(A) \subset \ww(A)$ due to Lemma~\ref{l:wsww}. Therefore,
$\ws(A) = \ww(A)$.
\end{proof}
Thanks to this theorem and the fact that
$\ww(X) \subset X$, we have
\[
\Aw=\ww(X) = \ws(X),
\]
provided $\Dc$ satisfies A1.
Note that if $\As$ exists, Theorem~\ref{l:wb=ab} implies that $\As = \ws(X)$.
Hence, we have the following.
\begin{corollary}
Let $\Dc$ be an evolutionary system satisfying A1. If $\As$ exists, then
\[
\As =\Aw.
\]
\end{corollary}

Similarly to the proof of Theorem~\ref{t:ws=ww,inv}, we can also obtain
the following extension of a
corresponding result for generalized semiflows (see \cite{B1}).
Note that due to the lack of concatenation, the proof of this theorem
strongly relies on property A1.
\begin{theorem}
Let $\Dc$ be an evolutionary system satisfying A1. Let $A\subset X$
be a weakly closed set. Then $A$ is invariant if and only if $A$
is positively invariant and quasi-invariant.
\end{theorem}
\begin{proof}
Clearly, if $A$ is positively invariant and quasi-invariant, then $A$ is invariant.
Assume now that $A$ is invariant. Then $A$ is positively invariant. To show that
it is quasi-invariant, consider any $x\in A$. Since $A$ is invariant, there exist 
$t_n \to \infty$ as $n\to \infty$ and $u_n \in \Dc([-t_n, \infty))$ with $u_n(-t_n) \in A$ and
$u_n(0) =x$. As in the proof of Theorem~\ref{t:ws=ww,inv}, passing to a subsequence
and dropping a subindex, we infer that there exists $u\in \Dc((-\infty,\infty))$ with
$u(0)=x$, such that $u_n \to u$ in $C([-T,\infty);\Xw)$ for all $T>0$. Since $A$ is
weakly closed, $u(t) \in A$ for all $t \in \mathbb{R}$, i.e., $A$ is quasi-invariant.
\end{proof}

Let
\[
\mathcal{I} :=\{ u_0: \ u_0=u(0) \mbox{ for some } u \in \Dc((-\infty, \infty))\}.
\]
Clearly, $\I$ is quasi-invariant and invariant. Moreover, it contains every
quasi-invariant an every invariant set.
Due to Theorem~\ref{t:ws=ww,inv}, it also follows that
\begin{equation} \label{e:invarinI}
\ww(A) \subset \I, \qquad \forall A \subset X.
\end{equation}

Now we will show that $\ww(A)$ captures a long-time behavior of every
trajectory starting in $A$, provided property A1 holds.

\begin{theorem}[Weak uniform tracking property] \label{t:longtimelim}
Let $\Dc$ be an evolutionary system satisfying $A1$. Let $A\subset X$.
Then for any $\epsilon >0$, there exists $t_0$, such that for any $t^*>t_0$,
every trajectory $u \in \Dc([0,\infty))$ with $u(0) \in A$ satisfies 
\[
\dd_{C([t^*,\infty);\Xw)} (u, v) < \epsilon,
\]
for some complete trajectory $v \in \Dc((-\infty,\infty))$ with $v(t) \in \ww(A)$
for all $t \in \mathbb{R}$.
\end{theorem}
\begin{proof}
Suppose the claim is not true. Then there exist $\epsilon>0$,
$u_n \in \Dc([0,\infty))$ with $u_n(0)\in A$, and
$t_n \to \infty$ as $n\to \infty$, such that
\begin{equation} \label{e:toobcontr}
d_{C([t_n,\infty);\Xw)}(u_n, v) \geq \epsilon,
\end{equation}
for all $n$, all $v \in \Dc((-\infty,\infty))$ with $v(t) \in \ww(A)$,
$t\in\mathbb{R}$.

On the other hand, consider a sequence
$v_n \in \Dc([-t_n,\infty))$,
$v_n(t)=u_n(t+t_n)$. Thanks to the fact that $\Dc([-t_n,\infty))$ is compact in
$C([-t_n,\infty); \Xw))$ for all $n$, using a diaganalization process,
passing to a subsequence and dropping a subindex, we infer that there exists
$v\in \Dc((-\infty, \infty))$, such that $v_n \to v$ in $C([-T,\infty); \Xw)$ as
$n \to \infty$ for all $T>0$.
In particular, $v(t) \in \ww(A)$ for all $t\in \mathbb{R}$. Finally, for large $n$
we have $d_{C([0,\infty);\Xw)}(v_n, v) < \epsilon$, which means that
$d_{C([t_n,\infty);\Xw)}(u_n, v(\cdot-t_n)) <  \epsilon$, a contradiction
with \eqref{e:toobcontr}.
\end{proof}

\begin{theorem}[Strong uniform tracking property] \label{t:strongtrack}
Let $\Dc$ be an asymptotically compact evolutionary system satisfying $A1$.
Let $A \subset X$.
Then for any $\epsilon >0$ and $T>0$, there exists $t_0$, such that for any
$t^*>t_0$, every trajectory
$u \in \Dc([0,\infty))$ with $u(0) \in A$ satisfies
\[
\ds(u(t), v(t)) < \epsilon, \qquad \forall t\in [t^*,t^*+T],
\]
for some complete trajectory $v \in \Dc((-\infty,\infty))$ with $v(t) \in \ws(A)$
for all $t \in \mathbb{R}$.
\end{theorem}
\begin{proof}
Suppose that the claim does not hold. Then there exist $\epsilon>0$, $T>0$,
and sequences $u_n \in \Dc([0,\infty))$ with $u_n(0) \in A$ and
$t_n \to \infty$ as $n\to \infty$, such that
\begin{equation} \label{e:stp1}
\sup_{t\in[t_n,t_n+T]}\ds(u_n(t),v(t)) \geq \epsilon, \qquad \forall n,
\end{equation}
for all $v \in \Dc((-\infty,\infty))$ with $v(t) \in \ws(A)$, $t\in \mathbb{R}$.

On the
other hand, Theorem~\ref{t:longtimelim} implies that there exists a sequence
$v_n \in \Dc((-\infty, \infty))$ with $v_n(t) \in \ww(A)$ for all $t$, such that
\begin{equation} \label{e:temdw}
\lim_{n\to \infty} \sup_{t\in[t_n,t_n+T]}\dw(u_n(t),v_n(t)) =0.
\end{equation}
Since $\Dc$ is asymptotically compact, Theorem~\ref{t:asymptoticcompact}
implies that $\ws(A)=\ww(A)$. Therefore, thanks to \eqref{e:stp1}, there exists a
sequence $\hat{t}_n \in [t_n,t_n+T]$, such that
\begin{equation} \label{e:temds}
\ds(u_n(\hat{t}_n),v_n(\hat{t}_n)) \geq \epsilon/2, \qquad \forall n.
\end{equation}
Due to the asymptotic compactness of $\Dc$, the sequences $\{u_n(\hat{t}_n)\}$
and $\{v_n(\hat{t}_n)\}$ are relatively strongly compact. Therefore, passing to subsequences and dropping a subindex, we obtain that $u_n(\hat{t}_n) \to x$, 
$v_n(\hat{t}_n) \to y$ strongly as $n\to \infty$ for some $x,y \in X$. Thanks to \eqref{e:temdw},
$x=y$, which contradicts \eqref{e:temds}.
\end{proof}

We conclude this section with a summary of the above results applied
to the weak global attractor $\Aw$. The tracking property of $\Aw$ connects
the global attractor with a trajectory attractor (see Section~\ref{s:trajectory}).

\begin{theorem} \label{c:weakA}
Let $\Dc$ be an evolutionary system satisfying A1. Then
$\Aw = \mathcal{I}$, and $\Aw$ is the maximal invariant and maximal
quasi-invariant set. Moreover,
for any $\epsilon >0$, there exists $t_0$, such that for any $t^*>t_0$,
every trajectory $u \in \Dc([0,\infty))$ satisfies 
\[
\dd_{C([t^*,\infty);\Xw)} (u, v) < \epsilon,
\]
for some complete trajectory $v \in \Dc((-\infty,\infty))$.

\end{theorem}
\begin{proof}

Note that $\Aw = \ww(X)$. Since $\I$ contains every
quasi-invariant set, Theorem~\ref{t:ws=ww,inv} implies that $\Aw \subset \I$.
On the other hand, thanks to \eqref{e:inclus}, it follows that
\[
\I \subset \ww(\I) \subset \ww(X) = \Aw.
\]
Therefore, $\Aw =\I$. The last statement of the theorem follows from
Theorem~\ref{t:longtimelim}.
\end{proof}

In the case $\Dc$ is asymptotically compact, this theorem together with
Theorems~\ref{t:asymptoticcompact} and \ref{t:strongtrack} implies the following.

\begin{theorem} \label{c:strongA1}
Let $\Dc$ be an asymptotically compact evolutionary system satisfying A1.
Then $\As = \mathcal{I}$ and $\As$ is the maximal invariant and maximal
quasi-invariant set. Moreover,
for any $\epsilon >0$ and $T>0$, there exists $t_0$, such that for any 
$t^*>t_0$, every trajectory $u \in \Dc([0,\infty))$ satisfies
\[
\ds(u(t), v(t)) < \epsilon, \qquad \forall t\in [t^*,t^*+T],
\]
for some complete trajectory $v \in \Dc((-\infty,\infty))$.
\end{theorem}

Note that the solutions to the 2D Navier-Stokes equations supplemented
with appropriate
boundary conditions form an asymptotically compact evolutionary system
satisfying A1. Therefore, Theorem~\ref{c:strongA1} yields the strong tracking
property for the 2D NSE.

\section{Evolutionary system with energy inequality}
In this section we will study an evolutionary system $\Dc$ satisfying
A2 and A3. In A2 we assume that $X$ is a bounded set in some uniformly
convex Banach space $H$ with
the norm denoted by $|\cdot|$, such that $\ds(x,y)=|x-y|$ for $x,y \in X$. We also
assume that for any $\epsilon >0$, there exists $\delta$, such that
for every $u \in \Dc([0,\infty))$ and $t>0$,
\[
|u(t)| \leq |u(t_0)| + \epsilon,
\]
for $t_0$ a.e. in $(t-\delta, t)$.

In A3, for every sequence $u_n \in \Dc([0,\infty))$, such that
$u_n \to u\in\Dc([0,\infty))$ in $C([0, T];\Xw)$ for some $T>0$, we assume that
$u_n(t) \to u(t)$ strongly a.e. in $[0,T]$.

 
\begin{theorem} \label{t:gaex}
Let $\Dc$ be an evolutionary system satisfying A2 and A3.
Let $u_n \in \Dc([T_1,\infty))$ be such that
$u_n \to  u$ in $C([T_1, T_2]; \Xw)$ as $n \to \infty$ for some $u\in\Dc([T_1,\infty))$.
If $u(t)$ is strongly continuous at some $t=t^* \in (T_1, T_2)$,
then $u_n(t^*) \to u(t^*)$ strongly in $X$.
\end{theorem}
\begin{proof}
Thanks to A3, there exists a set $E$ of measure zero, such that
$u_n(t) \to u(t)$ strongly on $[T_1,T_2]\setminus E$. Let $\epsilon >0$.
Due to the energy inequality $A2$ and strong continuity of $u(t)$,
there exists $t_0 \in [T_1,t^*)\setminus E$, such that
\[
|u_n(t^*)| \leq |u_n(t_0)| + \epsilon, \qquad |u(t_0)| \leq |u(t^*)| + \epsilon,
\]
for every $n$.
Taking the upper limit as $n\to \infty$, we obtain
\[
\begin{split}
\limsup_{n\to \infty} |u_n(t^*)| & \leq |u(t_0)| + \epsilon\\
& \leq |u(t^*)| + 2\epsilon.
\end{split}
\]
Since this inequality holds for an arbitrary $\epsilon$, we have
\[
\limsup_{n\to \infty} |u_n(t^*)| \leq |u(t^*)|.
\]
Hence, the weak convergence of $u_n(t^*)$ to $u(t^*)$ is in fact strong.
\end{proof}

\begin{theorem} \label{t:continmpacom}
Let $\Dc$ be an evolutionary system satisfying A1, A2, and A3.
If $\Dc((-\infty, \infty)) \subset C((-\infty,\infty); \Xs)$, then $\Dc$ is
asymptotically compact.
\end{theorem}
\begin{proof}
Take any sequences $t_n \to \infty$ and $x_n \in R(t_n)X$. 
Without loss of generality, there exists $T>0$, such that $t_n > T$ for all $n$.
Since $X$ is weakly compact, passing to a subsequence and
dropping a subindex, we can assume that $x_n \to x$ weakly as $n \to \infty$,
for some $x\in X$.

Since $x_n \in R(t_n)X$, there exist
$u_n \in \Dc([-t_n,\infty))$, such that $u_n(0) = x_n$.
Thanks to A1, using a diagonalization process,
passing to a subsequence and dropping a subindex,
we infer that there exists  $u\in \Dc((-\infty, \infty))$, such that
\[
u_n \to u \qquad \mbox{in} \qquad C([-T,T]; \Xw).
\]
Since $u(t)$ is
strongly continuous at $t=0$, Theorem~\ref{t:gaex} implies that
$x_n=u_n(0) \to x$ strongly as $n \to \infty$. Therefore, $\Dc$ is asymptotically
compact.
\end{proof}

This theorem, together with Theorem~\ref{t:asymptoticcompact},
immediately implies the following.

\begin{corollary} \label{c:strongA}
Let $\Dc$ be an evolutionary system satisfying A1, A2, and A3.
If every complete trajectory is strongly continuous, then $\Dc$ 
possesses a strongly compact strong global attractor $\As$.
\end{corollary}

It is clear from the above proof that for the strong convergence towards a weak omega
limit we only need the strong continuity of complete trajectories that pass through
the omega limit. More precisely, we have the following result.

\begin{theorem}
Let $\Dc$ be an evolutionary system satisfying A1, A2, and A3.
Let $A \subset X$ be such that there exists $u\in\Dc([0,\infty))$
with $u(0) \in A$.  Assume that $u(t)$ is strongly continuous at $t=0$
for every $u\in \Dc((-\infty,\infty))$ with $u(0) \in \ww(A)$.
Then $\ww(A)$ is a nonempty strongly compact set
that uniformly strongly attracts $A$. Moreover, $\ws(A)=\ww(A)$.
\end{theorem}
\begin{proof}
The proof follows the same lines as the proofs of
Theorems~\ref{t:acompomega} and \ref{t:continmpacom}.
\end{proof}

Now we show that the strong continuity of a limit trajectory implies
a uniform strong convergence towards the trajectory.
\begin{theorem} \label{t:unifstrongconv}
Let $\Dc$ be an evolutionary system satisfying A1, A2, and A3.
Let $u_n \in \Dc([0,\infty))$ be such that
$u_n \to u$ in $C([0, \infty); \Xw)$ as $n \to \infty$ for some 
$u \in \Dc([0, \infty))$. If
$u(t)$ is strongly continuous on $(0,\infty)$,
then $u_n \to u$ in $L^\infty_\mathrm{loc}((0,\infty);H)$.
\end{theorem}
\begin{proof}
Suppose that $u_n \nrightarrow u$ in $L^\infty_\mathrm{loc}((0,\infty);H)$.
Then passing to a subsequence and dropping a subindex we can assume that
there exist an interval $[a,b]\subset (0,\infty)$, $\epsilon >0$, and a sequence
$t_n \in [a,b]$, such that
\[
\ds(u_n(t_n), u(t_n)) \geq \epsilon, \qquad \forall n.
\]

Again passing to a subsequence and dropping a subindex, we obtain that
$t_n \to t_0 \in [a,b]$ as $n\to \infty$. Without loss of generality, assume
that $t_n \geq t_0/2$ for all $n$. Consider a sequence of trajectories
$v_n(t) = u_n(t-t_0+t_n)$, $t\in[t_0/2,\infty)$.
Note that $u_n \in \Dc([t_0/2,\infty))$ for all $n$.
Since $u(t)$ is continuous on $(0,\infty)$, it follows that
$\ds(u(t_n),u(t_0)) \leq \epsilon/2$ for $n$ large enough.
Therefore, since $v_n(t_0) = u_n(t_n)$, the triangle inequality implies that
\begin{equation} \label{e:cnucontr} 
\begin{split}
\ds(v_n(t_0), u(t_0))& \geq \ds(v_n(t_0), u(t_n)) - \ds(u(t_n), u(t_0))\\
& \geq \epsilon - \frac{\epsilon}{2} = \frac{\epsilon}{2},
\end{split}
\end{equation}
for $n$ large enough.

On the other hand, since
$u_n \to u$ in $C([0,\infty);\Xw)$ as $n\to \infty$, we have that
\[
\begin{split}
\dw(v_n(t),u(t)) &\leq \dw(v_n(t), u(t-t_0+t_n)) + \dw(u(t-t_0+t_n),u(t))\\
&=\dw(u_n(t-t_0+t_n), u(t-t_0+t_n)) + \dw(u(t-t_0+t_n),u(t))\\
&\to 0,
\end{split}
\]
as $n\to \infty$, for all $t \in [t_0/2, \infty)$. Therefore, thanks to A1, passing to a subsequence
and dropping a subindex, we obtain that $v_n \to u$ in $C([t_0/2, 2t_0];\Xw)$
as $n\to \infty$.
Since $u(t)$ is strongly continuous at $t=t_0$, Theorem~\ref{t:gaex}
implies that $v_n(t_0) \to u(t_0)$ strongly as $n\to \infty$, a contradiction
with \eqref{e:cnucontr}.
\end{proof}

Finally, we show that the strong continuity of a trajectory is equivalent
to the strong continuity from the right, provided A2 and A3 hold.

\begin{theorem} \label{t:gaex}
Let $\Dc$ be an evolutionary system satisfying A2 and A3.
Let $u \in \Dc([0,\infty))$ and $t^*>0$. Then $u(t)$ is strongly continuous
at $t=t^*$ if and only if $u(t)$ is strongly continuous from the left
at $t=t^*$.
\end{theorem}
\begin{proof}
Assume that $u(t)$ is strongly continuous from the left at $t=t^*$.
Let $\epsilon>0$. Due to the energy inequality A2, there exists $\delta>0$,
such that for every $t\in (t^*, t^*+\delta)$, there exists
a sequence $t_n <t^*$, $t_n \to t^*$ as $n \to \infty$, such that
\[
|u(t)| \leq |u(t_n)| + \epsilon,
\]
for every $n$. Thanks to the strong continuity from the left of $u(t)$ at
$t=t^*$, we obtain that
\[
|u(t)| \leq |u(t^*)| + \epsilon,
\]
for all $t\in(t^*,t^*+\delta)$.
Finally, taking the limit as $\epsilon \to \infty$, we obtain
\[
\limsup_{t \to t^*+} |u(t)| \leq |u(t^*)|.
\]
Hence, $u(t)$ is strongly continuous from the right at $t=t^*$. Therefore,
$u(t)$ is continuous at $t=t^*$.
\end{proof}

\section{Trajectory attractor} \label{s:trajectory}
A trajectory attractor for the 3D NSE was introduced in $\cite{S}$ and
further studied in \cite{CV, SY}.
In this section we define a trajectory attractor for the evolutionary
system $\Dc$ and discuss its properties that follow from the results in preceding
sections.

Consider an evolutionary system $\Dc$ satisfying A1.
Let $\F := C([0,\infty);\Xw)$.
As before, a function in $\Dc([0,\infty))$ is called a trajectory. Denote
\[
\Kp := \Dc([0,\infty)) \subset \F.
\]
Note that $\Kp$ is compact in $\F$ due to $A1$.
Define the translation operator $T(s)$
\[
(T(s)u)(t) := u(t+s)|_{[0,\infty)}, \qquad u \in \F.
\]
Due to the property 3 of the evolutionary system $\Dc$
(see Definition~\ref{Dc}), we have that
\[
T(s)\Kp \subset \Kp, \qquad \forall s\geq 0.
\]
For a set $P \subset \Kp$ define
\[
P(t) := \{u(t): u \in P\}, \qquad t \geq 0.
\]
Note that since we do not assume the uniqueness of the trajectories,
$P$ does not have
to contain all the trajectories starting at $P(0)$. More precisely,
by the definition of $R(t)$, for all $t\geq 0$ we have
\[
P(t)=(T(t)P)(0) \subset R(t)P(0), \qquad P \subset \Kp.
\]
On the other hand,
\begin{equation} \label{e:slice}
\Kp(t)=(T(t)\Kp)(0) = R(t)X, \qquad \forall t\geq0,
\end{equation}
since $\Kp$ includes all the trajectories in the evolutionary system.
For a set $P \subset \F$ and $r>0$ denote
\[
B(P,r) := \{u\in \F: \ d_{C([0,\infty);\Xw)}(u,P) < r\}.
\]
A set $P \subset \F$ uniformly attracts a set
$Q \subset \Kp$
if for any $\epsilon>0$ there exists $t_0$, such that
\[
T(t)Q \subset B(P, \epsilon), \qquad \forall t \geq t_0.
\]

\begin{definition}
A set $P \subset \F$ is a trajectory attracting set if it uniformly
attracts $\Kp$.
\end{definition}

\begin{lemma} \label{l:tatr-atr}
Let $P$ be a trajectory attracting set. Then $P(0)$ is a weakly attracting set.
\end{lemma}
\begin{proof}
Indeed, if $T(t)\Kp \subset B(P,\epsilon)$ for some $\epsilon$, $t\geq 0$,
then, thanks to \eqref{e:slice}, we have
\[
R(t)X =(T(t)\Kp)(0) \subset (B(P,\epsilon))(0) \subset \Bw(P(0),\epsilon).
\]
This concludes the proof.
\end{proof}

\begin{definition} A set
$\At \subset \F$ is a trajectory attractor if
$\At$ is a minimal compact trajectory attracting  set, and $T(t) \At = \At$ for all $t\geq 0$.
\end{definition}
It is easy to see that the intersection of two compact trajectory attracting sets
is a trajectory attracting set. Therefore, if a trajectory attractor exists, it
is unique.
Let $\K := \Dc((-\infty, \infty))$, which is called the kernel of $\Dc$. Let also
\[
\Pi_+ \K := \{u(\cdot)|_{[0,\infty)}:u \in \K\}.
\]

\begin{theorem}
Let $\Dc$ be an evolutionary system satisfying A1. Then the trajectory
attractor exists and
\[
\At = \Pi_+\K.
\]
\end{theorem}
\begin{proof}
Since $\Dc([-T,\infty))$ is compact in $C([-T,\infty);\Xw)$ for all $T>0$, using a diagonalization process, we obtain that $\Pi_+\K$ is compact in $\F$.
Moreover, due to Theorem~\ref{c:weakA}, $\Pi_+\K$ uniformly attracts $\Kp$.

Now assume that there exits a compact trajectory attracting set $P$ strictly
included in $\Pi_+\K$. Then there exist $\epsilon>0$ and
\[
u \in  \Pi_+\K \setminus B(P,\epsilon).
\]
Let $v\in \Dc((-\infty,\infty))$ be such that $v|_{[0,\infty)}=u$. Let also
$v_n(\cdot)=v(\cdot-n)|_{[0,\infty)}$. Note that $v_n \in \Kp$ and
\[
T(n)v_n = u \notin B(P,\epsilon), \qquad \forall n.
\]
Therefore, $P$ is not a trajectory attracting set, a contradiction.

Finally, the properties 2 and 4 of $\Dc$ immediately imply that
$T(t) \Pi_+\K = \Pi_+\K$ for all $t\geq 0$.
\end{proof}

Note that Theorem~\ref{c:weakA} also yields that
\[
\Aw = \At(t), \qquad \forall t \in \mathbb{R}.
\]
Finally, we will show that the strong continuity of complete trajectories implies
a uniform strong convergence of solutions toward the trajectory attractor.

\begin{theorem}
Let $\Dc$ be an evolutionary system satisfying A1, A2, and A3.
If $\At \subset C([0,\infty);\Xs)$, then the trajectory attractor $\At$
uniformly attracts $\Kp$ in $L^\infty_\mathrm{loc}((0,\infty); H)$. 
\end{theorem}
\begin{proof}
Since $\At \subset C([0,\infty);\Xs)$, Theorem~\ref{t:continmpacom} implies
that the evolutionary system $\Dc$ is asymptotically compact. Therefore,
Theorem~\ref{c:strongA1} yields that $\At$ uniformly attracts $\Kp$ in
$L^\infty_\mathrm{loc}((0,\infty); H)$.

\end{proof}


\section{3D Navier-Stokes equation}
In this section we will apply the above results to the space periodic 3D incompressible
Navier-Stokes equations (NSE)
\begin{equation} \label{NSE1}
\left\{
\begin{aligned}
&\ddt u - \nu \Delta u + (u \cdot \nabla)u + \nabla p = f,\\
&\nabla \cdot u =0,\\
&u, p, f \mbox{ are periodic with period } L \mbox{ in each space variable,}\\
&u, f \mbox{ are in } L^2_{\mathrm{loc}}(\mathbb{R}^3)^3,\\
\end{aligned}
\right.
\end{equation}
where $u$, the velocity, and $p$, the pressure, are unknowns; $f$ is
a given driving force, and $\nu>0$ is the kinematic  viscosity coefficient
of the fluid. By a Galilean change of
variables, we can assume that the space average of $u$ is
zero, i.e.,
\[
\int_\Omega u(x,t) \, dx =0, \qquad \forall t,
\]
where $\Omega=[0,L]^3$ is a periodic box.

First, let us introduce some notations and functional setting.
Denote by $(\cdot,\cdot)$ and $|\cdot|$ the $L^2(\Omega)^3$-inner product and the
corresponding $L^2(\Omega)^3$-norm.
Let $\mathcal{V}$ be the space of all $\mathbb{R}^3$ trigonometric polynomials of
period $L$ in each variable satisfying 
$\nabla \cdot u =0$ and $\int_\Omega u(x) \, dx =0$.
Let $H$ and $V$ to be the
closures of $\mathcal{V}$ in $L^2(\Omega)^3$ and $H^1(\Omega)^3$, respectively.
Define the strong and weak distances by
\[
\ds(u,v):=|u-v|, \qquad
\dw(u,v)= \sum_{\kappa \in \mathbb{Z}^3} \frac{1}{2^{|\kappa|}}
\frac{|u_{\kappa}-v_{\kappa}|}{1 + |u_{\kappa}-v_{\kappa}|},
\qquad u,v \in H,
\]
where $u_{\kappa}$ and $v_{\kappa}$ are Fourier coefficients of $u$
and $v$ respectively.

Let also  $P_{\sigma} : L^2(\Omega)^3 \to H$ be the $L^2$-orthogonal
projection, referred to as the Leray projector. Denote by
$A=-P_{\sigma}\Delta = -\Delta$ the Stokes operator with the domain
$D(A)=(H^2(\Omega))^3 \cap V$. The Stokes operator is a self-adjoint
positive operator with a compact support.
Let
\[
\|u\| := |A^{1/2} u|,
\]
which is called the enstrophy norm.
Note that $\|u\|$ is equivalent to the $H^1$-norm of $u$ for $u\in D(A^{1/2})$.

Now denote $B(u,v):=P_{\sigma}(u \cdot \nabla v)\in V'$ for all 
$u, v \in V$. This bilinear form has the following property:
\[
(B(u,v),w)=-(B(u,w),v), \qquad u,v,w \in V,
\]
in particular, $(B(u,v),v)=0$ for all $u,v \in V$.

Now we can rewrite (\ref{NSE1}) as the following differential equation in $V'$:

\begin{equation} \label{NSE}
\ddt u + \nu A u +B(u,v) = g,\\
\end{equation}
where $u$ is a $V$-valued function of time and $g = P_{\sigma} f$. Throughout, we will assume that
 $g$ is time independent and $g\in H$.

\begin{definition}
A weak solution  of  \eqref{NSE1} on $[T,\infty)$ (or $(-\infty, \infty)$, if
$T=-\infty$) is an $H$-valued
function $u(t)$ defined for $t \in [T, \infty)$, such that
\[
\ddt u \in L_{\mathrm{loc}}^1([T, \infty); V'), \qquad
u(\cdot) \in C([T, \infty); \Hw) \cap 
L_{\mathrm{loc}}^2([T, \infty); V),
\]
and
\[
\left(\ddt u, v\right) + \nu ((u, v)) + (B(u,u),v) =(g, v)
\qquad \mbox{a.e. in }t,
\forall v \in V.
\]
\end{definition}

\begin{theorem}[Leray, Hopf] \label{thm:Leray}
For every $u_0 \in H$, there exists a weak solution of (\ref{NSE1}) on $[T,\infty)$ with $u(T)=u_0$
satisfying the following energy inequality
\begin{equation} \label{EI}
|u(t)|^2 + 2\nu \int_{t_0}^t \|u(s)\|^2 \, ds \leq
|u(t_0)|^2 + 2\int_{t_0}^t (g(s), u(s)) \, ds
\end{equation}
for all $t \geq t_0$, $t_0$ a.e. in $[T,\infty)$.
\end{theorem}

\begin{definition} \label{d:ex}
A Leray-Hopf solution of \eqref{NSE1} on the interval $[T, \infty)$
is a weak solution on $[T,\infty)$ satisfying the
energy inequality (\ref{EI}) for all $T \leq t_0 \leq t$,
$t_0$ a.e. in $[T,\infty)$. The set $Ex$ of measure $0$ on which the energy
inequality does not hold will be called the exceptional set.
\end{definition}




It is known that there exists an
absorbing ball for the 3D Navier-Stokes equations (see, e.g., \cite{ConF}).

\begin{proposition} \label{p:aset}
The 3D Navier-Stokes equations possess an absorbing ball
\[
B = B_{\mathrm{s}}(0, R),
\]
where $R$ as any number larger that $|g| \nu^{-1} L/(2\pi)$.
\end{proposition}

Let $X$ be a closed absorbing ball
\[
X= \{u\in H: |u| \leq R\},
\]
which is also weakly compact. Then for any bounded set $A \subset H$,
there exists a time $t_0$, such that
\[
u(t) \in X, \qquad \forall t\geq t_0,
\]
for every Leray-Hopf solution $u(t)$ with the initial data $u(0) \in A$.
Classical NSE estimates (see \cite{ConF}) imply that for any sequence of Leray--Hopf solutions
$u_n$ (not only for the ones guaranteed by Theorem~\ref{thm:Leray})
the following result holds.

\begin{lemma} \label{l:convergenceofLH}
Let $u_n(t)$ be a sequence of Leray-Hopf solutions of  \eqref{NSE1},
such that $u_n(t) \in X$ for all $t\geq t_0$. Then 
\[
\begin{aligned}
u_n \ \ &\mbox{is bounded in} \ \ L^2([t_0,T];V),\\
\ddt u_n \ \  &\mbox{is bounded in} \ \ L^{4/3}([t_0,T];V'),
\end{aligned}
\]
for all $T>t_0$.
Moreover, there exists a subsequence $u_{n_j}$ of $u_n$ that converges
in $C([t_0, T]; \Hw)$ to some Leray-Hopf solution $u(t)$, i.e.,
\[
(u_{n_j},v) \to (u,v) \qquad
\mbox{uniformly on} \qquad  [t_0,T],
\]
as $n_j\to \infty$, for all $v \in H$.
\end{lemma}

Consider an evolutionary system for which
a family of trajectories consists
of all Leray-Hopf solutions of the 3D Navier-Stokes equations
in $X$. More precisely, define
\[
\begin{split}
\Dc([T,\infty)) := \{&u(\cdot): u(\cdot)
\mbox{ is a Leray-Hopf}
\mbox{ solution on } [T,\infty)\\
&\mbox{and } u(t) \in X \ \forall t \in [T,\infty)\},
\qquad T \in \mathbb{R},
\end{split}
\]
\[
\begin{split}
\Dc((\infty,\infty)) := \{&u(\cdot): u(\cdot)
\mbox{ is a Leray-Hopf} \ \mbox{ solution on } (-\infty,\infty)\\
&\mbox{and } u(t) \in X \ \forall t \in (-\infty,\infty)\}.
\end{split}
\]

Clearly, the properties 1--4 of $\Dc$ hold.
Therefore, thanks to Corollary~\ref{thm:exofA}, the weak global attractor $\Aw$
for this evolutionary system exists. Moreover, we have the following.
\begin{lemma} \label{l:compact}
The evolutionary system $\Dc$ of the 3D NSE
satisfies A1, A2, and A3.
\end{lemma}
\begin{proof}
First note that $\Dc([0,\infty)) \subset C([0,\infty);\Hw)$ by the definition of
a Leray-Hopf solution. Now take any sequence
$u_n \in \Dc([0,\infty))$, $n=1,2, \dots$.
Thanks to Lemma~\ref{l:convergenceofLH}, there exists
a subsequence, still denoted by  $u_n$, that converges
to some $u^{1} \in \Dc([0,\infty))$ in $C([0, 1];\Hw)$ as $n \to \infty$.
Passing to a subsequence and dropping a subindex once more, we obtain that
$u_n \to u^2$ in $C([0, 2];\Hw)$ as $n \to \infty$ for some 
$u^{2} \in \Dc([0,\infty))$.
Note that $u^1(t)=u^2(t)$ on $[0, 1]$.
Continuing
this diagonalization process, we obtain a subsequence $u_{n_j}$
of $u_n$ that converges
to some $u \in \Dc([0,\infty))$ in $C([0, \infty);\Hw)$ as $n_j \to \infty$.
Therefore, A1 holds.

Now, given $\epsilon>0$, let $\delta=\epsilon/(2 |g|R)$. Take any $u \in \Dc([0,\infty))$
and $t>0$. Since $u(t)$ is a Leray-Hopf solution,
it satisfies the energy inequality \eqref{EI}
\[
|u(t)|^2 + 2\nu \int_{t_0}^t \|u(s)\|^2 \, ds \leq
|u(t_0)|^2 + 2\int_{t_0}^t (g, u(s)) \, ds,
\]
for all $0 \leq t_0 \leq t$, $t_0 \in [0,\infty) \setminus Ex$, where
$Ex$ is a set of zero measure. Hence,
\[
\begin{split}
|u(t)|^2 &\leq |u(t_0)|^2 + 2(t-t_0) |g| R\\
&\leq |u(t_0)|^2 + \epsilon,
\end{split}
\]
for all $t_0\geq 0$, such that $t_0 \in (t-\delta,t)\setminus Ex$. Therefore, A2 holds.

Let now $u_n \in \Dc([0,\infty))$ be such that $u_n \to u\in\Dc([0,\infty))$ in
$C([0, T];\Xw)$ as $n\to \infty$ for some
$T>0$. Thanks to Lemma~\ref{l:convergenceofLH}, the sequence $\{u_n\}$
is bounded in $L^2([0,T];V)$. Hence,
\[
\int_{0}^T |u_n(s)-u(s)|^2 \, ds \to 0, \qquad \mbox{as}
\qquad  n \to \infty.
\]
In particular, $|u_n(t)| \to |u(t)|$ as $n \to \infty$ a.e. on $[0,T]$,
i.e., A3 holds.
\end{proof}

Now Theorem~\ref{c:weakA} and Corollary~\ref{c:strongA} yield the following.
\begin{theorem} The weak global attractor $\Aw$ for the 3D
NSE exists, $\Aw$ is the maximal invariant set, and
\[
\Aw = \ww(X)=\ws(X)=\{u(0): u \in \Dc((-\infty, \infty))\}.
\]
Moreover, for any $\epsilon >0$ there exists $t_0$,
such that for any $t^*>t_0$,
every Leray-Hopf solution $u \in \Dc([0,\infty))$ satisfies
\[
\dd_{C([t^*,\infty);\Xw)} (u, v) < \epsilon,
\]
for some complete trajectory $v \in \Dc((-\infty,\infty))$.

\end{theorem}

\begin{theorem}
If every complete trajectory of the 3D NSE
is strongly continuous, then the weak global
attractor is a strongly compact strong global attractor. 
In addition, for any $\epsilon >0$ and $T>0$, there exists $t_0$,
such that for any $t^*>t_0$,  every Leray-Hopf solution $u \in \Dc([0,\infty))$
satisfies
\[
\ds(u(t), v(t)) < \epsilon, \qquad \forall t\in [t^*,t^*+T],
\]
for some complete trajectory $v \in \Dc((-\infty,\infty))$.
\end{theorem}

Finally, we note that all the other results from the previous sections apply to the 3D
Navier-Stokes equations as well.


\begin{thebibliography}{99}

\bibitem{BV}
A.V. Babin and M.I. Vishik, Maximal attractors of semigroups corresponding to evolution 
differential equations. {\it Mat. Sb.} {\bf 126} (1985), 397Ð419; English transl: {\it Math. USSR 
Sb.} {\bf 54} (1986), 387-408.

\bibitem{B1} J. M. Ball, Continuity properties and global attractors of generalized
semiflows and the Navier-Stokes equations, {\it J. Nonlinear Sci.} {\bf7} (1997),
475Ð502. Erratum: {\it J. Nonlinear Sci.} {\bf 8} (1998),  233.

\bibitem{B2} J. M. Ball, Global attractors for damped semilinear wave equations,
{\it Discr. Cont. Dyn. Sys.} {\bf10} (2004), 31Ð52.

\bibitem{Bar} E. A. Barbashin, On the theory of generalized dynamical systems,
{\it Moskov. Gos. Ped. Inst.  U\v{c}en. Zap.}, {\bf 2} (1948), 110Ð133.
English translation by U.S. Department of Commerce, Office 
of Technical Services, Washington D.C. 20235.

\bibitem{CMR} T. Caraballo, P. Mar'n-Rubio, and J. C. Robinson,
A comparison between two theories for multi-valued semiflows and their
asymptotic behaviour, {\it Set-Valued Anal.} {\bf 11} (2003), 297--322.

\bibitem{CV} V. V. Chepyzhov and M. I. Vishik, {\it Attractors for Equations of
Mathematical Physics}, American Mathematical Society Colloquium Publications
{\bf 49}, American Mathematical Society, Providence, RI, 2002.

\bibitem{C} A. Cheskidov, Blow-up in finite time for the dyadic model of the
Navier-Stokes equations, {\it Trans. Amer. Math. Soc.}, to appear, arXiv:math.AP/0601074.

\bibitem{CF} A. Cheskidov and C. Foias, On global attractors of the 3D
Navier-Stokes equations, {\it J. Diff. eq.}, to appear, arXiv:math.AP/0608475.

\bibitem{CFP1} A. Cheskidov, S. Friedlander, and N. Pavlovi\'{c},
A dyadic model for the inviscid fluid equations: the global attractor,
preprint.


\bibitem{ConF} P. Constantin and C. Foias, {\it Navier-Stokes Equation}, University of Chicago Press, Chicago, 1989. 

\bibitem{FMRT} C. Foias, O. P. Manley, R. Rosa, and R. Temam, {\it Navier-Stokes
equatinon and Turbulence}, Encyclopedia of Mathematics and its Applications
{\bf 83}, Cambridge University Press, Cambridge, 2001. 

\bibitem{FT85}  C. Foias and R. Temam, The connection between the Navier-Stokes
equations, and turbulence theory, {\it Directions in Partial
Differential Equations} (Madison, WI, 1985), Publ. Math. Res. Center
Univ. Wisconsin, 55-73.

\bibitem{FP}
S. Friedlander and N. Pavlovi\'c, Blowup in a three-dimensional vector
model for the Euler equations, {\it Comm. Pure Appl. Math.} {\bf 57} (2004),
705--725.


\bibitem{H} J. K. Hale, {\it Asymptotic behavior of dissipative systems}, 
Amer. Math. Soc., Providence, RI, 1988.

\bibitem{HLS}J. K. Hale,  J. P. LaSalle, and M. Slemrod,
Theory of a general class of dissipative processes, {\it J. Math. Anal. Appl.}
{\bf 39} (1972), 177--191.

\bibitem{KP}
N. H. Katz and N. Pavlovi\'c, Finite time blow-up for a dyadic model of the Euler
equations, {\it Trans. Amer. Math. Soc.} {\bf 357} (2005), 695--708.


\bibitem{L} O. Ladyzhenskaya, {\it Attractors for semigroups and evolution equations},
Cambridge University Press, Cambridge, 1991.

\bibitem{MV} V. S. Melnik and J. Valero,  On attractors of multivalued semi-flows
and differential inclusions, {\it Set-Valued Anal.} {\bf 6} (1998), 83--111.

\bibitem{LR}  J.A. Langa and J. C. Robinson, Determining Asymptotic Behaviour from the Dynamics on Attracting Sets, {\it Journal of Dynamics and Differential Equations} {\bf 11} (1999), 319-331.

\bibitem{O} L. Onsager, Statistical Hydrodynamics, {Nuovo Cimento (Supplemento)} {\bf{6}} (1949),  279--287.

\bibitem{R} R. M. S. Rosa,  Asymptotic regularity condition for the strong
convergence towards weak limit sets and weak attractors of the 3D
Navier-Stokes equations, {\it J. Diff. eq.}, to appear.

\bibitem{S} G. R. Sell, Global attractors for the three-dimensional Navier-Stokes
equations, {\it J. Dynam. Differential Equations} {\bf 8} (1996), 1-33. 

\bibitem{SY} G. R. Sell and Y. You, {\it Dynamics of evolutionary equations},
Applied Mathematical Sciences {\bf 143}, Springer-Verlag, New York, 2002.


\bibitem{T2} R. Temam, {\it Infinite Dimensional Dynamical Systems in Mechanics and
Physics}, Applied Mathematical Sciences {\bf 68}, (2nd Edition, 1997) Springer Verlag, New York, 1988.

\end{thebibliography}
\end{document}